\documentclass[11pt]{amsart}
\usepackage{amsfonts,amsmath,amsthm,amscd,amssymb,latexsym,cite,verbatim,texdraw,floatflt,
caption2}
\RequirePackage{hyperref}
\usepackage[a4paper]{geometry}

\title{ On a question of Clark and Ledet  }
\author{ Igor  Protasov}

\address{I.Protasov: Taras Shevchenko National University of Kyiv, Department of Computer Science and Cybernetics, Academic Glushkov pr. 4d, 03680 Kyiv, Ukraine}
\email{i.v.protasov@gmail.com}

\begin{document}
\begin{abstract} 
Given a $T$-sequence on a countable abelian group $G$,
we prove that there exists $2^{2^{|G|}}$ Hausdorff group topologies in which this sequence converges to $0$. 
This answers a question posed in Intern. J. Math. Math. Sci. {\bf 24}(3) (2000), 145-148.
\end{abstract}
\maketitle



\vspace{5 mm}

All groups under consideration are supposed to be abelian and all topologies are supposed to be Hausdorff. 

A sequence $(a_n)_{n\in \omega}$ in a group $G$ is called a {\it $T$-sequence} if there exists a group topology on $G$ in which $(a_n)_{n\in \omega}$ converges to $0$. 
In \cite{b2}, the authors called a sequence in $G$ to be  a 
{\it $T_\Omega$-sequence} if there exists an uncountable  set of group topologies in which this sequence converges to $0$, and asked does there exist a 
$T$-sequence on an infinite group which is not a $T_\Omega$-sequence? We give  the negative answer after short explanations. 

\vspace{6 mm}

For a group topology $\tau$ on $G$, $\mathcal{N}_\tau$ denotes the family of all neighborhoods of $0$ in
 $\tau$. 
Two non-discrete group topologies $\tau_1$, $\tau_2$ on $G$ are called {\it transversal} (or complemented) if 
the upper bound $\tau_1 \vee \tau_2$ of $\tau_1$ and 
$\tau_2$ is the discrete topology. 
Equivalently,  $\tau_1$ and $\tau_2$ are transversal if 
there exist $U_1 \in \mathcal{N}_{\tau_1}$ and $U_2 \in \mathcal{N}_{\tau_2}$ such that 
$U_1 \cap U_2 = \{0\}$.
If $\tau_1$, $\tau_2$ are transversal then the lower bound $\tau_1 \wedge \tau_2$ is Hausdorff  and the family $\{ U + V : U\in \mathcal{N}_{\tau_1}$, $V\in \mathcal{N}_{\tau_2}\}$ is a base for 
$\mathcal{N}_{\tau_1 \wedge \tau_2}$.

\vspace{5 mm}
We use the following observation.

\vspace{5 mm}

$(\ast) \ $  {\it Let $\tau$, $\tau^\prime$, $\tau_1$, $\tau_2$ be group topologies on $G$ such that $\tau$,
$\tau^\prime$ are transversal, 
$\tau^\prime \subseteq \tau_1$, 
$\tau^\prime \subseteq \tau_2$. 
If $\tau\wedge  \tau_1 \subseteq \tau\wedge  \tau_2$
then $\tau_1  \subseteq \tau_2$.
}

\vspace{5 mm}

To verify $(\ast)$, we pick $U\in \mathcal{N}_\tau $,
$ \ V\in \mathcal{N}_\tau^\prime$, 
$ \ V_0 \in \mathcal{N}_{\tau^\prime} \ $
such that 
$U\cap V=\{0\}$,
$ \ V_0 = - V_0$,  $ \ V_0 + V_0 \subseteq  V$.

Let $V_1\in \mathcal{N}_{\tau_1}$,  
$V_1\subseteq V_0 $.
Since $\tau\wedge\tau_1 \subseteq \tau\wedge\tau_2$,
there exist $V_2\in \mathcal{N}_{\tau_2}$
such that $V_2\subseteq V_0 $ and $V_2\subseteq U+V_1 $.
For any $v_2\in V_2$, we take $u\in U$, $v_1 \in V_1$
such that $v_2 = u+ v_1$.
Since $U\cap V = \{ 0\}$ and $V_0 + V_0 \subseteq V$,
we have $v_1 = v_2$ so $V_2 \subseteq V_1$ and $\tau_1\subseteq \tau_2$.

\vspace{5 mm}

{\bf Theorem.} {\it For every $T$-sequence $(a_n)_{n\in \omega}$ in an infinite group $G$, there exists $2^{2^{|G|}}$ group topologies in which $(a_n)_{n\in \omega}$ converges to $0$.
\vspace{5 mm}

Proof.} We denote by $\tau$ the strongest group topology in which $(a_n)_{n\in \omega}$ converges to $0$. If all but finitely many $a_n$
are equal to $0$ then $\tau$ is discrete and we apply 
\cite{b1}, so we suppose that $\tau$ is non-discrete and  consider two cases. 

\vspace{4 mm}

{\bf Case 1:} $G$ is countable. By [4, Corollary 2.4.10], there exists a group topology $\tau^\prime$
transversal to $\tau$. Since $G$ is countable, after some weakening, we may suppose that $\mathcal{N}_{\tau^\prime}$ has a countable base. Using arguments from the proof of Theorem 1.1 in 
\cite{b3}, we can choose an injective sequence $(b_n)_{n\in\omega}$ which converge to $0$ and such that,
for every partition $\omega=W_1 \cup W_2$ into infinite subsets $W_1$, $W_2$, the strongest group topologies in which $(b_n)_{n\in W_1}$,  $(b_n)_{n\in W_2}$ converge to $0$ are transversal. 
Then, for every  free ultrafilter $\varphi$  on $\omega$, we denote by $\tau_\varphi$ the supremum of the group topologies determined by the sequences 
$\{ (b_n)_{n\in F} : F \in \varphi \}$.
If $\varphi , \ \psi$ are distinct ultrafilters then 
$\tau_\varphi , \ \tau_\psi$ are transwersal. 
By $(\ast)$, 
$\tau \wedge \tau_\varphi  \ \neq \  \tau \wedge\tau_\psi$.
Since $(a_n)_{n\in \omega}$ converges to $0$, $\tau$ in every topology $\tau \wedge\tau_\varphi$
and the number of free ultrafilter on $\omega$ is $2^ {\mathfrak{c}}$, this case is proven.

\vspace{4 mm}

{\bf Case 2:} $G$ is uncountable. We denote by $H$ the subgroup generated by the set $\{ a_n : n\in \omega \}$.
Since $H$ is countable and $G$ is a subgroup of the direct sum of countable groups (the rationals and quasi-cyclic),  we can choose a subgroup $S$ of $G$
such that $|S|=|G|$
and $S\cap H = \{ 0 \}$.
We denote by $M$ the set of all group topologies on $G$
in which $H$ is open. 
By \cite{b1}, $|M|= 2^{2^{|H|}}$ and it suffices to 
note that $(a_n)_{n\in\omega}$ converges to $0$
in each topology $\tau\wedge \mu$, $\mu\in M$.
$ \  \Box$

\end{document}